\documentclass[12pt,a4paper]{article}
\usepackage{makeidx}
\usepackage{graphicx}
\usepackage{latexsym}
\usepackage{amsmath}
\usepackage{amsfonts}
\usepackage{verbatim}
\usepackage{nameref}
\usepackage{hyperref}
\usepackage{sectsty}
\usepackage{enumerate}

\usepackage{amssymb}

\usepackage{amsmath}
 \usepackage[dvips]{epsfig}
\usepackage{amsfonts}
\usepackage{amsmath,amscd}
\usepackage{amsbsy}
\usepackage{amscd}
\usepackage{float}
\usepackage{rotating}
\usepackage{xy} \xyoption{all} \xyoption{import} \xyoption{rotate}
\usepackage{pst-plot}
\usepackage{psfrag}
\usepackage{a4wide}
\newcommand{\bX}{{\bf x}}
\newcommand{\WW}{\hat{W}}

\newcommand{\alp}{\alpha}
\newcommand{\VV}{V}
\newcommand{\Oo}{\Omega}
\newcommand{\dO}{{\rm d} \Omega}
\newcommand{\dG}{{\rm d} \Gamma}
\newcommand{\Gu}{{\Gamma_u}}
\newcommand{\Gt}{{\Gamma_t}}
\newcommand{\calB}{{\cal B}}
\newcommand{\uu}{{u}}
\newcommand{\baruu}{\bar{u}}
\newcommand{\be}{{\bf e}}
\newcommand{\cc}{{c}}
\newcommand{\xx}{{x}}
\newcommand{\la}{\label}
\newcommand{\half}{\frac{1}{2}}
\newcommand{\bt}{{t}}
\newcommand{\bB}{{\bf B}}
\newcommand{\bE}{{\bf E}}
\newcommand{\bI}{{\bf I}}
\newcommand{\bC}{{\bf C}}
\newcommand{\bQ}{{\bf Q}}
\newcommand{\bF}{{\bf F}}
\newcommand{\bT}{{\bf T}}
\newcommand{\bn}{{\bf n}}
\newcommand{\bx}{{\bf x}}

\newcommand{\sta}{{\rm sta}}

\newcommand{\bxi}{\mbox{\boldmath$\xi$}}
\newcommand{\bchi}{\mbox{\boldmath$\chi$}}
\newcommand{\bvsig}{\zeta}
\newcommand{\barzeta}{\bar{\zeta}}
\newcommand{\bveps}{\xi}
\newcommand{\bsig}{\mbox{\boldmath$\tau$}}
\newcommand{\bgamma}{\mbox{\boldmath$\gamma$}}
\newcommand{\btau}{{\mbox{\boldmath$\tau$}}}

\newcommand{\calP}{{\cal P}}

\newcommand{\calE}{{\cal{E}}}
\newcommand{\calS}{{\cal{S}}}
\newcommand{\calU}{{\cal{U}}}
\newcommand{\calT}{{\cal{T}}}
\newcommand{\calX}{{\cal{U}}}

\newcommand{\real}{\mathbb{R}}
\newcommand{\eb}{\begin{equation}}
 \newtheorem{thm}{Theorem}
\newcommand{\ee}{\end{equation}}

 \newtheorem{remark}{Remark}

\def \div{\mbox{div\hskip 1pt}}

\def \tr{\mbox{tr\hskip 1pt}}

\def\bm#1{\mbox{\boldmath{$#1$}}}   

\title{\bf  Analytical  Solutions to General  Anti-Plane Shear Problems In Finite Elasticity}

\author{David Yang  Gao  \\[0.2cm]
\small   Federation University Australia, Mt Helen, VIC 3353, Australia  \\
\small  Research School of Engineering, Australian National University, Canberra, Australia\\
 [0.2cm]
}

\date{}

\begin{document}
\maketitle
\begin{center}
{\em Dedicated to Professor Ray Ogden on the occasion of his 70th birthday}
\end{center}
\begin{abstract}
 This paper presents  a pure complementary energy variational method for  solving  anti-plane  shear
problem in finite elasticity. Based on the canonical duality-triality  theory developed by the author,
  the
nonlinear/nonconex partial differential equation for the large  deformation problem is converted into an
algebraic equation in dual space, which can, in
principle, be solved   to obtain  a complete set of stress solutions.
Therefore, a general analytical solution form of the deformation is obtained 
subjected to a compatibility condition.
Applications are illustrated by
 examples with both convex and  nonconvex stored strain energies governed by
 quadratic-exponential and power-law material models, respectively.
Results show  that  the nonconvex variational problem could have multiple solutions at each
 material point, the complementary  gap function and the triality theory can be
used to identify both global and local extremal solutions, while the popular
(poly-, quasi-, and rank-one) convexities  provide only local minimal criteria,
 the   Legendre-Hadamard condition does not guarantee uniqueness of solutions.
This paper demonstrates again  that the pure
complementary energy principle and the  triality  theory play important
  roles in finite deformation theory and nonconvex analysis.
\end{abstract}

\section{Introduction}
Anti-plane shear deformation problems arise naturally from many real world applications,
such as rectilinear steady flow of simple fluids \cite{fos-ser}, interface stress effects of nanostructured materials \cite{luo-wang}, structures with cracks \cite{psm}, layered/composite functioning materials \cite{narita-shindo,yu-yang},   and phase transitions in solids \cite{silling}.
During the past half century, such problems in finite deformation theory
have been subjected to extensively  study
by both mathematicians and engineering scientists \cite{ball,gurtin-temam,hill,hill2,knowles,pucci,pucci-sacc}.
 As  indicated in the review article by C.O. Horgan \cite{horgan}),
 anti-plane shear deformations are one of the simplest classes of deformations that solids can undergo.
 In anti-plane shear (or longitudinal shear, generalized shear) of a cylindrical body, the displacement is parallel to the generators of the cylinder and is independent of the axial coordinate.
 In recent years, considerable attention has been paid to the analysis of anti-plane shear deformations within the context of various constitutive theories (linear and nonlinear) of solid mechanics. Such studies were largely motivated by the promise of relative analytic simplicity compared with plane problems since the governing equations are a single second-order  partial differential equation rather than higher-order or coupled systems of partial differential equations. Thus the anti-plane shear problem plays a useful role as a pilot problem, within which various aspects of solutions in solid mechanics may be examined in a particularly simple setting.

Generally speaking, the anti-plane shear problem in linear elasticity is  governed by linear partial differential equation, which can be solved easily by well-developed  analytical  methods.
However,    in finite elasticity, the associated variational problem is usually nonlinear or nonconvex,
 which could have multiple solutions. Traditional methods for solving nonconvex variational problems are proved
 to be very difficult, or even impossible.
 The well-known generalized convexities  and  {\em Legendre-Hadamard condition} can be used only for identifying local minimal solutions.
 Numerical methods (such as FEM and FDM, etc) for solving nonconvex variational problems lead to a global optimization problem.
 Due to the lacking  of  global optimality condition, most of nonconvex optimization problems are considered to be NP-hard
 in nonconvex optimization and computer science \cite{gao-sherali-amma,murty}. Extensive research has been
 focused on solving such nonconvex optimization problems,  and a special  research field, i.e., the global optimization
 has been developed during the past 15 years \cite{gao-sherali}.

Complementary variational principles and methods play important roles in continuum mechanics.
It is  known that in finite deformation theory, the  Hellinger-Reissner principle (see \cite{hell-14, reiss53})
 and
the Fraeijs de Veubeke principle (see \cite{veub72}) hold for both convex and nonconvex problems.
 But,  these well-known principles are not considered as the {\em pure
complementary variational principles} since the Hellinger-Reissner principle
involves both the displacement field and the second Piola-Kirchhoff stress tensor,
and  the Fraeijs de Veubeke principle
needs both the rotation tensor and the first Piola-Kirchhoff
stress as its variational arguments.
Therefore,  the question about the existence of a pure complementary variational
principle in general finite deformation theory was argued
  for several decades
(see \cite{koiter76, lee-shie80, lee-shie80b, oden-redd83,ogden75,ogden77}).
A systematic study   on the   invariant
conditions for various complementary energy functionals in elasticity
was  given by Li and Gupta \cite{li-cupta}. 
Also, since the extremality condition
is a fundamentally difficult problem in nonconvex variational analysis and global optimization,
all the  classical
 complementary-dual variational principles and associated numerical methods
can't be used for solving nonconvex variational/optimization problems  in finite deformation theory.

Canonical duality-triality  is a newly developed and powerful methodological theory,
which is composed mainly of  a {\em canonical   transformation}, a pure complementary-dual energy variational principle,  and a
{\em triality theory}.
The canonical   transformation can be used to model complex systems within a unified framework and to establish
perfect dual problems in nonconvex analysis and global optimization.
 The pure complementary-dual variational principle shows that a class of nonlinear partial differential equations
 are equivalent to  certain  algebraic equations  which can be solved to obtain analytical solutions in stress space.
 The triality theory comprises  a {\em canonical min-max duality} and a pair of {\em double-min, double-max dualities}.
 The canonical min-max duality  can be used to identify
  global minimizer;  while the double-min and double-max dualities can be used to identify  local minimizer  and local maximizer, respectively.
 The original idea of the canonical duality theory
 was introduced by Gao and Strang in general nonconvex variational problems
 \cite{gao-strang89a}. The triality theory was discovered in post-buckling analysis of a large deformed beam model \cite{gao-amr97}.
The pure complementary principle was first proposed in 1999 \cite{gao-mrc99}, which  has been used successfully for solving
  finite deformation problems.
  In a set of papers published recently by Gao and Ogden \cite{gao-ogden-qjmam,gao-ogden-zamp},
  it is shown that by using this theory, complete sets of analytical solutions can be obtained
  for one-dimensional nonlinear/nonconvex problems. Their results illustrated an  important fact
  that smooth analytic
or numerical solutions of a nonlinear mixed boundary-value problem might not be minimizers
of the associated potential variational problem.
  For global optimization problems in finite dimensional space,
  the canonical duality theory has been used successfully for
and global optimization,
see \cite{gao-sherali-amma}

The purpose of this paper is to illustrate the application of this
  pure complementary variational principle and the triality theory by solving nonlinear and
nonconvex variational problems in anti-plane shear deformation. The
remainder article is organized as the following. The next section
discusses the finite anti-plane  shear deformation and constitutive
laws. Based on the equilibrium equation and a general constitutive law, a
nonlinear potential variational problem is formulated.
 Section 3 shows how this nonlinear potential variational
problem can be transformed as a canonical dual problem such that  a pure complementary energy principle can be obtained,
and by which, how the nonlinear partial differential equation for deformation
can be converted into an algebraic equation in stress space, so that an analytical solution form for the displacement
can be formulated.
This section also shows how the global optimal solution can be identified by the triality theory.
Section 4 presents an application to convex problem governed by a  quadratic-exponential stored energy,
which has a unique solution.
 While for
nonconvex strain energy, Section 5 shows that the boundary value
problem is not equivalent to the variational problem. By using the
canonical dual transformation and the pure complementary
variational principle, the nonlinear differential equation can be
converted into a cubic algebraic equation, which possesses at most
three real roots. Therefore a complete set of solutions to the
potential variational problem is obtained. The triality theory can
be used to identify global and local minimizers.

  The results presented in this paper show  that  the pure complementary
  energy principle and the triality theory
  are potentially useful  in finite deformation theory.

\section{Anti-plane shear deformation and variational problem\label{azimdef}}

Consider a homogeneous, isotropic elastic cylinder $\calB \subset \real^3$
with generators parallel to the $\be_3$ axis and
with cross section a sufficiently nice region $\Oo \subset \real^2$ in the
$\be_1 \times \be_2$ plane.
The so-called anti-plane shear  deformation is defined by
\begin{equation}\label{defor}
\bchi = \bX + \uu(\xx_1, \xx_2) \be_3, \;\; \forall (\xx_1, \xx_2 ) \in \Oo
\end{equation}
where $ (\xx_1, \xx_2, \xx_3 )$ are cylindrical  coordinates in the
reference configuration (assumed free of stress) $\calB$ relative to a
cylindrical  basis $\{\be_\alp\},\alp = 1,2,3$, and $\uu:\Oo \rightarrow \real$ is
 the amount of shear (locally a simple shear) in the planes normal to $\be_3$.

  We suppose that on boundary $\Gu \subset \partial \Oo$ the homogenous boundary condition is given, i.e.,
\begin{equation}\label{27}
\uu(\bx)= 0 \;\; \forall \bx  \in \Gu,
\end{equation}
while on the boundary $\Gt = \partial \Oo  \cap \Gu$, the shear force   is prescribed
\[
{\bf t} (\bx) =  t(\bx)  \be_3 \;\; \forall \bx \in \Gt.
\]

The deformation gradient tensor, denoted $\mathbf{F}$,
can be readily calculated for the deformation of the form (\ref{defor}):
\begin{equation}\label{1}
\mathbf{F}= \nabla \bchi = \mathbf{I}+ \be_3  \otimes (\nabla \uu) = \left(\begin{array}{ccc}
1 & 0 & 0 \\
0 & 1& 0 \\
\uu_{,1} & \uu_{,2} & 1 \end{array} \right),
\end{equation}
where  $\uu_{,\alpha}$ represents $\partial \uu/\partial x_\alpha$ for $\alpha = 1,2$,
$\mathbf{I}$ is the identity tensor; while the corresponding left Cauchy-Green tensor, denoted
$\mathbf{B}$, is
\begin{equation}\label{25}
\mathbf{B}=\mathbf{FF}^{\rm
T}=
\left( \begin{array}{ccc}
1 & 0 & \uu_{,1}\\
0 & 1 & \uu_{,2} \\
\uu_{,1} & \uu_{,2} & 1 + |\nabla \uu |^2
\end{array} \right),
\end{equation}
where     the notation  $^{\rm T}$ indicates the transpose (of a
second-order tensor).

The principal invariants of $\mathbf{B}$, denoted $I_1,I_2,I_3$,
are defined by
\begin{equation}
I_1=\tr\mathbf{B},\quad  I_2 = \half [ (\tr \bB)^2 - \tr(\bB)^2 ]
 ,\quad I_3=\det\mathbf{B},
\end{equation}
and, for the considered anti-plane shear problem, these reduce to
\begin{equation}
I_1= I_2 = 3 + | \nabla \uu |^2 , \;\;
 \quad I_3=1 .\la{invariant}
\end{equation}
In this paper,  the notation $|\nabla \uu|^2 = \uu_{,1}^2 + \uu_{,2}^2$ represents the Euclidean norm in $\real^2$.



In view of \eqref{invariant}, if we let $\bgamma = \nabla u$,  we may now introduce a new function,
denoted $\hat{W}$, such that the stored-energy density function $W(\bF)$ can be written as
\begin{equation}\label{22}
W(\bF) =  \hat{W}(\bgamma)
\end{equation}
and the   associated in-plane Cauchy stress tensor $\bsig$
\begin{equation}
\bsig  =\frac{\partial
\hat{W}}{\partial\bgamma}\label{shear-stress}
\end{equation}
is  the shear stress.

For the deformation and constitutive law discussed in the
foregoing sections the equilibrium equation
$\div\bm\bsig = {0}$ (in the absence of body forces) has
just one  non-trivial components, namely the equilibrium  equation
\begin{equation}\label{11}
\nabla \cdot \btau =
\frac{ \partial \tau_{1}}{ \partial
\xx_1}+ \frac{ \partial \tau_{2}}{ \partial
\xx_2}=0
\end{equation}

Let   the kinetically admissible
 space be  defined by
\eb
\calX_a = \{ \uu(\bx)  \in \hspace*{1pt} {\cal C} [{\bar{\Oo}}; \real] \;\; | \;\; \nabla \uu  \in {\cal L}^{2p} [\bar{\Oo}; \real^2],  \;\;
\uu(\bx) = 0  \; \forall \bx \in \Gu\},
\ee
where $\bar{\Oo} = \Oo \cup \partial \Oo$ represents the
closure of the set $\Oo$,  ${\cal L}^{2p}$ is the standard notation of Lebesgue integrable space with $p\in\real$, which will be discussed in Section 5.
Then the  minimal potential energy principle leads to
 the following variational problem (primal problem $(\calP)$ for short) for the
determination of the deformation function $\uu$:

\begin{equation}
(\calP): \;\;  \min\left\{
\Pi (\uu)=   \int_\Oo \hat{W}(\nabla \uu ) \dO - \int_\Gt  \bt \uu  \dG  \;\;\; |\;\; \uu \hspace*{1pt}\in
 \calX_a     \right\} , \label{eq-vp1}
\end{equation}
where $\min \{ * \}$ represents for finding minimum value of the statement in $\{ * \}$.

The criticality condition $\delta \Pi(\uu) = 0 $ leads to the
following mixed boundary value problem:
\begin{equation}
(BVP)  : \;\; \left\{
\begin{array}{l}
\nabla \cdot  \frac{\partial \hat{W} (\nabla \uu)}{\partial (\nabla \uu)}   = 0   \;\;   \forall \bx \in  \Oo, \vspace{.2cm}\\
 \bn\cdot \frac{\partial \hat{W} (\nabla \uu)}{\partial (\nabla \uu)}   = \bt    \;   \; \forall \bx \in \Gt .
\end{array} \right.
\end{equation}
Since the stored energy $\hat{W}(\bgamma)$ is a nonlinear
 function of the shear strain, the $(BVP)$ may possess multiple solutions and
  each solution represents a stationary point
 of the total potential energy $\Pi(\uu)$. Therefore, if $\Pi(\uu)$ is nonconvex, the
  variational problem $(\calP)$ is not equivalent the boundary value problem $(BVP)$.
Traditional direct methods for solving the nonlinear  boundary value problem are usually very difficult, also
 how to identify the global minimizer in nonconvex analysis is a fundamentally difficult task.
 It turns out that in general nonlinear elasticity,
even some qualitative questions such as regularity and stability are considered as outstanding open problems
 \cite{ball}.

 On the other hand, let the statically admissible space to be
 \eb
 \calT_a = \{ \btau \in {\cal C} [\bar{\Oo}; \real^2] | \; \nabla \cdot \btau(\bx) = 0  \;\; \forall \; \bx \in \Oo, \;\;
 \bn \cdot \bsig(\bx)   = \bt  \;  \forall \bx \in  \Gt \}.
\ee
The
   complementary variational problem can be
stated as the following:
\eb
\min \left\{ \Pi^*(\bsig  ) =  \int_\Oo\hat{W}^*(\bsig )
  {\rm d} \Oo  \;\; | \;\; {\btau \in \calT_a } \right\},
  \ee
where $\hat{W}^*(\btau)$ is the so-called complementary energy
function (or density), defined by the Legendre transformation:
\eb
\hat{W}^*(\bsig ) = \{ \bgamma \cdot  \bsig - \hat{W}(\bgamma) \; | \;\;
\bsig  = \nabla \hat{W}(\bgamma)\}.
\ee

In finite deformation theory, if the  strain energy density $\hat{W}(\bgamma)$ is  nonconvex,
the Legendre conjugate $\hat{W}^*(\bsig)$ can not be uniquely obtained \cite{ogden,sewell}.
In this case, the classical complementary energy variational principle can not be used for
solving nonconvex finite deformation problems.
Although by the Fenchel transformation
\eb
\WW^\sharp(\btau) = \sup \{  \bgamma \cdot  \bsig - \hat{W}(\bgamma) \; | \;\; \bgamma \in {\cal L}^{2p}[\bar{\Oo}; \real^2]\},
\ee
the Fenchel conjugate $\WW^\sharp(\btau)$ is always convex, and the Fenchel-Moreau dual problem can be obtained as
\eb
\max \left \{ \Pi^\sharp = - \int_\Oo \WW^\sharp(\btau) \dO \;\; | \; \; \btau \in \calT_a \right \},
\ee
 the Fenchel-Young inequality
\[
\WW(\bgamma) \ge \bgamma \cdot  \bsig - \hat{W}^\sharp (\btau ) \;\;
\]
leads to
\eb
\theta = \min_{\uu \in \calX_a}  \Pi(\uu) -  \max_{\btau \in \calT_a } \Pi^\sharp(\btau) \ge 0
\ee
and the non zero $\theta \neq 0 $ is  the well-known duality gap in nonconvex analysis.
According to Sir M. Atiyah \cite{atiyah}, duality in mathematics is not a theorem, but a ``principle".
Therefore, the duality gap is not allowed in mathematical physics. It turns out that
 the existence of a pure complementary energy principle in finite elasticity
was a well-known open problem which has been discussed for over thirty years.
This problem was solved in 1999 \cite{gao-mrc99} when a complementary energy principle was proposed in
terms of the first and second Piola-Kirchhoff stresses only.

In the following  sections, we will demonstrate the
application of the canonical duality theory
and the pure complementary variational principle for solving the proposed variational
problem.
 In order to examine this problem in
detail, we will consider the energy function in both convex and nonconvex forms.

\section{Canonical dual problem and extremality theory}
The key step of the canonical dual transformation
  is to introduce a new
 geometrical measure $\bveps = \Lambda (\uu)$ and a
 canonical function $\VV(\bveps)$
   such that the stored energy function
   $\WW(\nabla \uu ) = \VV(\Lambda(\uu))$.
 By the definition introduced in \cite{gao-dual00} that a real-valued function
  $\VV(\bveps)$ is called a canonical function if
  the duality relation
  \eb
  \bvsig = \frac{\partial \VV(\bveps)}{\partial \bveps}
  \ee
 is invertible such that the conjugate function $\VV^*(\bvsig)$ of $\VV(\bvsig)$
 can be defined uniquely
 by the Legendre transformation:
 \eb
 \VV^*(\bvsig) = \left\{ \bveps   \bvsig - \VV(\bveps) \; | \; \bvsig = \frac{\partial \VV(\bveps)}{\partial \bveps}\right\}
 .
 \ee

 The canonical dual transformation has a solid foundation in physics.
 According to the {\em frame-invariance axiom} \cite{ciarlet,holz,ogden},
 instead of the linear deformation   $\gamma = \nabla \uu$,
 the strain  energy  $\hat{W}(\gamma)$ should be a  function of a quadratic measure $\xi = \Lambda(\uu)$.
In view of (\ref{invariant}) and (\ref{22}),  for this anti-plane shear deformation problem
 we can simply choose the geometrical measure
 $\xi = \Lambda(\uu) = \half |\nabla \uu |^2$, which is a quadratic mapping
 from $\calX_a$ to a closed convex set
 \eb
 \calE_a =
 \{ \xi \in {\cal L}^{p}[\bar{\Oo}; \real]  |\; \xi(\bx)  \ge 0  \;\; \forall \bx \in \Oo\}.
\ee

  Let $\calE^*_a$ be the range of the canonical duality mapping $\nabla \VV:\calE_a \rightarrow \calE^*$,
  such that on $\calE_a \times \calE^*_a$, the following
  canonical duality relations hold:
  \eb
  \bvsig = \nabla \VV(\bveps) \;\; \Leftrightarrow \;\; \xi = \nabla \VV^*(\bvsig) \;\;
  \Leftrightarrow \;\; \VV(\xi) + \VV^*(\bvsig) = \xi \bvsig. \label{eq-cdr}
  \ee

In the terminology  of finite elasticity, if the geometrical measure $\xi$ can be viewed as a
Cauchy-Green type strain,  its conjugate
 $\zeta $ is a second Piola-Kirchhoff type stress.
 For many hyper-elastic materials, the stored energy function   could be nonconvex in
 the deformation gradient, but is usually convex in the Cauchy-Green type strain measure.
 Thus, replacing  $\WW(\nabla \uu) $ in the
 total potential energy $\Pi(\uu)$ by
 the canonical form $\VV(\Lambda(\uu)) $, the primal problem $(\calP) $ can be written in
 the following canonical form:
 \eb
 (\calP): \;\; \min \left\{ \Pi(\uu) = \int_\Oo \VV(\Lambda(\uu)) \dO - \int_\Gt \uu \bt \dG | \; \; \uu \in \calX_a \right\}.
 \ee
 Furthermore, by using the Fenchel-Young equality
 $\VV(\Lambda(\uu))
 = \Lambda(\uu)   \bvsig - \VV^*(\bvsig)$,
   the so-called {\em total complementary energy functional} originally  proposed  by Gao and Strang
   in \cite{gao-strang89a} can be written  as
 \eb
 \Xi(\uu, \bvsig) = \int_{\Omega} [\half |\nabla \uu|^2  \bvsig - \VV^*(\bvsig)
   ] {\rm d} \Omega - \int_{\Gamma_t} \uu   \bt {\rm d} \Gamma. \label{xig}
 \ee
This two-field functional is well-defined on $\calX_a \times \calE^*_a$.
Let
\eb
\calS^+ = \{ \zeta \in \calE^*_a | \; \; \zeta(\bx) \ge 0\;\; \forall \bx \in \Oo \}.
\ee
The following theorem is a special case of the general result by Gao and Strang \cite{gao-strang89a}.
  \begin{thm}[Complementary-Dual Variational Extremum  Principle]
 If $(\baruu, \barzeta)$ is a critical point of $\Xi(\uu, \zeta)$,
 then $\baruu$ is a local solution to $(BVP)$.
 Moreover, if $\VV(\xi)$ is convex and $\barzeta  \in \calS^+  $, then
 $\baruu$ is a global optimal solution to the minimal variational problem $(\calP)$ and
 \eb
 \Pi(\baruu) = \min_{\uu \in \calX_a } \Pi(\uu) =  \Xi(\baruu, \barzeta) =\min_{\uu \in \calX_a} \max_{\zeta \in \calS^+}
 \Xi(\uu, \zeta) = \max_{\zeta \in \calS^+} \min_{\uu \in \calX_a}
 \Xi(\uu, \zeta) . \label{minmaxp}
 \ee
 \end{thm}


 {\bf Proof}.  By the criticality condition $\delta \Xi(\baruu, \barzeta) = 0$, we obtain
 \eb
 \Lambda(\baruu) = \half |\nabla \baruu |^2 = \nabla \VV^*(\barzeta), \;\;\mbox{ in } \Oo \label{eq-inconstitutive}
 \ee
 \eb
 \nabla \cdot (\barzeta \nabla \baruu) = 0 \;\; \mbox{ in } \Oo, \;\; \bn \cdot (\barzeta \nabla \baruu) = t
 \; \mbox{ on } \Gt. \label{eq-balance}
 \ee
 By the canonical duality (\ref{eq-cdr}), we know that the Euler-Lagrangian equation (\ref{eq-inconstitutive}) is
 equivalent to the canonical constitutive-geometrical equation $\barzeta = \nabla \VV(\Lambda(\baruu))$.
 Combining this with the equilibrium equation (\ref{eq-balance}), we know that the critical point $(\baruu, \barzeta)$
 solves the boundary value problem $(BVP)$ and $\baruu$ is a critical point of the the total potential
 energy  $\Pi(\uu)$.

 By the convexity of the canonical energy $\VV(\xi)$, we have (see \cite{gao-strang89a})
 \[
 \VV(\xi) - \VV(\bar{\xi}) \ge (\xi - \bar{\xi} ) \nabla \VV(\bar{\xi}) \;\; \forall \xi, \bar{\xi} \in \calE_a.
 \]
 Let $\xi = \Lambda(\uu)$ ,  $\bar{\xi} = \Lambda(\baruu)$, and $\barzeta = \nabla \VV(\Lambda(\baruu))$,  we obtained
 \[
\Pi(\uu) - \Pi(\baruu) \ge \int_\Oo [  \barzeta ( \Lambda(\uu) - \Lambda(\baruu)) ]\dO - \int_\Gt t (\uu - \baruu) \dG
\;\; \forall \uu \in \calX_a.
\]
Let  $\uu = \baruu + \delta \uu$. By the fact that $\Lambda(\uu)$ is a quadratic operator, we have (see \cite{gao-strang89a})
\[
\Lambda(\uu) =  \Lambda(\baruu + \delta \uu)   = \Lambda(\baruu) + (\nabla  \delta \uu )^T (\nabla \baruu) + \Lambda(\delta \uu).
\]
Therefore, if $(\baruu, \barzeta)$ is a critical point of $\Xi(\uu, \zeta)$ and $\barzeta \in \calS^+$, we have
\[
\Pi(\uu) - \Pi(\baruu) = G(\delta \uu, \barzeta) =  \int_\Oo \barzeta \Lambda(\delta \uu) \dO \ge 0 \;\; \forall \delta \uu.
\]
This shows that $\baruu$ is a global minimizer of $\Pi(\uu)$ over $\calX_a$. \hfill $\Box$\\

\begin{remark}[Gao-Strang's Gap Function  and Global Optimality Condition]   $\;$ \newline
{\em   Theorem 1 is actually the direct application of the general result of  \cite{gao-strang89a},
 and
 \[
 G(\uu, \zeta) = \int_\Oo \zeta \Lambda(  \uu) \dO
 \]
 is   the so-called  complementary gap function first introduced by Gao and Strang in 1989 \cite{gao-strang89a}.
 Since the geometrical operator $\Lambda(\uu) = \half |\nabla \uu|^2$ is  quadratic,
the gap function
\[
G(\uu, \zeta ) \ge 0 \; \;\forall \uu \in\calU_a   \;\; \mbox{  if and only if } \;\; \zeta \in \calS^+ .
\]
 Therefore,   the total complementary energy
 $\Xi(\uu, \zeta)$ is a saddle functional on  $\calX_a \times \calS^+$, i.e.
 \[
 \Xi(\uu, \barzeta) \ge \Xi(\baruu, \barzeta) \ge \Xi(\baruu, \zeta) \;\; \forall (\uu, \zeta)
 \in \calX_a \times \calS^+.
 \]
 Thus, by the canonical min-max duality, we have (\ref{minmaxp}).
 Theorem 1 shows that the gap function $G(\uu, \barzeta) \ge 0 $ provides a global optimality condition
 for the nonconvex variational problem $(\calP)$. This gap function also plays a key role in global optimization
 (see \cite{gao-sherali-amma}). Based on this complementary extremum principle and the general
 canonical primal-dual  mixed finite element method \cite{cai-gao-qin,gao-yu},
  an efficient algorithm can be developed for solving general anti-plane shear problems.   }
  \end{remark}

By the virtual work principle, for any given statically admissible  $\btau \in \calT_a$, we have
\eb
\int_\Oo (\nabla \uu) \cdot \btau \dO = \int_{\Gamma} ( \btau \cdot \bn)  \uu \dG -
 \int_\Oo (\nabla \cdot \btau) \uu \dO  = \int_\Gt  \bt \uu \dG  \;\;
\; \forall  \uu  \in \calX_a  . \label{virtual}
\ee
Replacing  the boundary integral in (\ref{xig}) by (\ref{virtual}),
 the total complementary energy functional $\Xi(\uu, \zeta)$ can be written
as
\eb
\Xi_\btau (\uu, \zeta) = \int_\Oo \left [ \half |\nabla \uu |^2 \zeta -  \VV^*(\zeta) - (\nabla \uu) \cdot \btau \right] \dO.
\ee
\begin{thm}
For any given statically admissible $\btau \in \calT_a$,
if $(\baruu, \barzeta)$ is a critical point of $\Xi_\btau(\uu, \zeta)$, then it is also a critical point of $\Xi(\uu, \zeta)$,  and
\eb
\Pi(\baruu) = \Xi(\baruu, \barzeta) = \Xi_\btau (\baruu, \barzeta) \;\; \forall \btau \in \calT_a . \label{xi=xit}
\ee
\end{thm}

{\bf Proof}.
For a given $\btau \in \calT_a$, the criticality condition $\delta \Xi_\btau(\baruu, \barzeta) = 0$ gives to
the inverse constitutive law $
\half |\nabla \baruu|^2 = \nabla \VV^*(\barzeta)   \;\;
$
and the balance equations
\eb
\nabla \cdot (\barzeta \nabla \baruu )  =  \nabla \cdot \btau \; \mbox{ in } \Oo,
\;\; \bn \cdot (\barzeta \nabla \baruu )  =  \bn \cdot \btau \;\; \mbox{ on } \Gt. \label{tuz}
\ee
Therefore,  $(\baruu, \barzeta)$ is a critical point of $\Xi (\uu, \zeta)$, and also a solution to $(BVP)$
for any given  $\btau \in \calT_a$. The equality (\ref{xi=xit}) can be proved easily by the virtual work principle and the canonical duality relations
(\ref{eq-cdr}).
\hfill $\Box$\\

 This theorem  shows that the statically admissible stress field $\btau \in \calT_a$ does not
change the  value of the functional $\Xi_\btau(\uu, \zeta)$. Note that from the criticality conditions
(\ref{tuz}), we have
\[
\zeta \nabla \uu = \btau,
\]
which shows the relation between the canonical stress and the first Piola-Kirchhoff stress.
By substituting  $\nabla \uu = \btau /\zeta$ in $\Xi_\btau(\uu, \zeta)$,
 the pure complementary energy
$\Pi^d(\zeta)$ can be obtained by the canonical dual transformation
  \cite{gao-mecc99}
\begin{equation}
\Pi^d(\zeta)  =  \sta \left\{   \Xi_\btau (\uu, \zeta) | \;  \forall {\uu \in \calX_a }\right\}
 =   - \int_\Oo \left(\frac{ \;\; |\bsig|^2   }{2   \zeta} + \VV^{* }
 (\zeta) \right) \dO,
\end{equation}
which is well-defined on
\eb
\calS^+_a = \{ \zeta \in \calS^+ | \;\; |\btau|^2/\zeta   \in {\cal L}[\bar{\Oo}; \real]  \}.
\ee
Therefore, the complementary variational problem which  is canonically dual to the potential variational problem
$(\calP)$  can be proposed as
  \eb
  (\calP^d): \;\; \max \left\{   \Pi^d(\zeta) =  - \int_\Oo \left(\frac{ \;\; |\bsig|^2   }{2   \zeta} + \VV^{* }
 (\zeta) \right) \dO   | \;\;  \zeta \in \calS^+_a \right\}.
  \ee
  According to \cite{gao-dual00}, we have the following result.
 \begin{thm}[Pure Complementary Energy Principle]
For a given statically admissible $\btau \in \calT_a$, if $(\baruu, \barzeta)$ is a critical point of
$\Xi_\btau(\uu, \zeta)$, then
 $\barzeta$ is a critical point of  $\Pi^d(\zeta)$, $\baruu$ is a critical point of $\Pi(\uu)$, and
 \eb
 \Pi(\baruu) = \Xi_\btau(\baruu, \barzeta) = \Pi^d(\barzeta).
 \ee
If $\VV(\xi)$ is convex, then
 $\baruu$ is a  global minimum solution to $(\calP)$ if and only if $\barzeta \in \calS^+_a$ is
 a solution to $(\calP^d)$, i.e.
 \eb
 \Pi(\baruu) = \min_{\uu \in \calX_a} \Pi(\uu) \;\; \Leftrightarrow \;\; \max_{\zeta \in \calS^+_a} \Pi^d(\zeta) = \Pi^d(\barzeta) . \label{thm-minmax}
 \ee
 The problem $(\calP)$ has a unique solution if $\barzeta (\bx) > 0, \; \forall \bx \in \Oo$.
 \end{thm}

{\bf Proof}. By the equation (\ref{xi=xit}) we know that for a given $\btau \in \calT_a$,
the functionals $\Xi(\uu, \zeta)$ and $\Xi_\btau(\uu, \zeta)$ have the same critical points set.
Particularly,  the criticality condition condition
 $\delta \Pi^d(\barzeta) = 0$ leads to
 \eb
 \frac{ |\btau|^2 }{2 ( \barzeta )^2} = \nabla \VV^*(\barzeta), \; \label{eq-deuler}
 \ee
 which is in fact the inverse constitutive-geometrical equation  (\ref{eq-inconstitutive}) subject to
 \eb
 \nabla \baruu =  \frac{\btau}{\barzeta}. \label{eq-u=ts}
 \ee
 Since $\btau$ is statically admissible, therefore,   $\barzeta \nabla \baruu = \btau$ satisfies
 equilibrium conditions (\ref{eq-balance}).
 This proved that  the critical point $\barzeta$ of the canonical dual problem $(\calP^d)$
 and the associated $\baruu$ are also critical point of $\Xi$.

Again by the canonical duality  (\ref{eq-cdr}),
we have
\[
\Xi(\baruu, \barzeta) = \int_\Oo   \VV(\Lambda(\baruu)) \dO - \int_\Gt t \baruu \dG
=\int_\Oo \WW(\nabla \baruu) \dO  - \int_\Gt t \baruu \dG = \Pi(\baruu).
\]
Dually, by using (\ref{eq-u=ts}), we have for any given $\btau \in \calT_a$
\[
\Xi(\baruu, \barzeta) = \Xi_\btau(\baruu, \barzeta) = \Pi^d(\barzeta).
\]
By the fact that $\Xi_\btau(\uu, \zeta)$ is a saddle functional on $\calX_a \times \calS^+_a$,
we have
\[
\min_{\uu \in \calX_a} \Pi(\uu) = \min_{\uu \in \calX_a} \max_{\zeta \in \calS^+_a} \Xi_\btau(\uu,\zeta) =
 \max_{\zeta \in \calS^+_a} \min_{\uu \in \calX_a}\Xi_\btau(\uu,\zeta) = \max_{\zeta \in \calS^+_a} \Pi^d(\zeta).
 \]
Thus, $ \baruu \in \calX_a $   is a  global minimum solution to $(\calP)$ if and only if
$ \barzeta \in \calS^+_a$  is a solution to the canonical dual problem $(\calP^d)$.
Moreover, if $\barzeta(\bx) > 0 \;\; \forall \bx \in \Oo$, then the gap function $G(\uu, \barzeta) > 0 \;\;
\forall \uu \neq 0$. From the proof of Theorem 1 we know that
the total potential energy $\Pi(\uu)$ is strictly convex on $\calX_a$, therefore, the global min is unique.
\hfill $\Box$\\

\begin{remark}
{\em  This theorem  shows that the complementary energy variational problem
$(\calP^d)$ is canonically dual to the potential variational problem $(\calP)$, i.e.
there is no duality gap.
The   canonical dual Euler-Lagrangian equation (\ref{eq-deuler}) shows that
the criticality condition  of the pure complementary energy
is an algebraic equation
\eb
\barzeta^2 \nabla \VV^*(\barzeta ) = \half |\btau|^2, \label{eq-cde}
\ee
which can be solved easily for many real applications. Therefore, the pure complementary energy principle plays
an important role in stress analysis and design.
But, for each $\btau$, the solution $\barzeta$
can only produce the deformation gradient
  $\nabla \baruu = \barzeta^{-1} \btau$.
  In order to obtain the primal solution $\baruu$ by solving the canonical dual problem,
  additional compatibility condition is needed.
  }
  \end{remark}

  \begin{thm}[Analytical  Solution Form]\label{thm-anas}
For a given statically admissible stress  $\btau (\bx)  \in \calT_a$ such that
 $\barzeta(\bx) $ is a solution of  the canonical dual equation (\ref{eq-cde}),
 then the vector-valued function
 \eb
  \nabla \baruu = \barzeta^{-1} (\bx) \btau
  \ee
  is a deformation solution to
 the $(BVP)$.

 Moreover,  if $\btau \in\calT_a $ is a potential field  and
\eb
 \btau \times ( \nabla \barzeta )  = 0\;\;\; \forall \bx \in \Oo, \label{eq-compatible}
 \ee
  then
        the path-independent line integral
 \eb
  \bar{\uu}(\bx) = \int_{\bx_0}^\bx  \bar{\zeta}^{-1}    \btau  \cdot  {\rm d} \bx \;\; \forall \bx_0 \in \Gu \label{eq-solu1}
  \ee
   is a solution of the boundary value problem $(BVP)$.

  \end{thm}

  {\bf Proof}. First,  by using chain role we know that  if $\nabla \baruu$ is a solution to $(BVP)$, it must satisfy
\[
\frac{\partial \hat{W}(\nabla \baruu)}{\partial (\nabla \baruu)} = \frac{\partial \VV(\Lambda(\baruu))}{\partial \Lambda(\baruu)}
\frac{\partial \Lambda(\baruu)}{\partial \nabla \baruu} = \barzeta \nabla \baruu \in \calT_a .
\]
Let $\btau =\barzeta \nabla \baruu  \in \calT_a$, we have
\eb
 \nabla \baruu = \barzeta^{-1} \btau ,
\ee
which is indeed the critical condition $\delta_{\uu} \Xi_\btau (\baruu, \barzeta) = 0$.
By the canonical duality (\ref{eq-cdr}),
 we have
  \[
  \half |\nabla \baruu |^2 = \half \frac{|\btau|^2 }{\barzeta^2} = \nabla \VV^*(\barzeta),
  \]
  which is  the canonical dual algebraic equation, i.e., the criticality condition of $\delta \Pi^d(\barzeta) = 0$.
  Therefore, for a given $\btau \in \calT_a$, if $\barzeta$ is  a solution to this canonical dual equation (\ref{eq-cde}),  then $ \nabla \baruu = \barzeta^{-1} \btau $ is the deformation field of the
   $(BVP)$.

Moreover, if $\baruu$ can be solved by the path-independent line integral (\ref{eq-solu1}),
the integrant $ \barzeta^{-1} \btau $ must be a potential field on $\Oo$, i.e.
$\nabla \times ( \barzeta^{-1} \btau ) = 0$.
This leads to
\[
\barzeta (\nabla \times  \btau)  + \btau \times (\nabla \barzeta) = 0 \;\; \mbox{ on } \Oo.
\]
Since $\btau(\bx)$ is a potential field on $\Oo$, i.e., there exists a scale-valued function
$\phi(\bx)$ such that $\btau = \nabla \phi(\bx)$, then we have
$\nabla \times \btau = \nabla \times (\nabla \phi) \equiv 0 $ on $\Oo$.
Therefore, as long as the condition $\btau \times ( \nabla \barzeta )  = 0 $ holds on $\Oo$,
the deformation gradient  $ \nabla \baruu = \barzeta^{-1} \btau $ is a potential field
and the displacement $\baruu$ can be obtained by the path-integral (\ref{eq-solu1}).
By the fact that   $\baruu(\bx_0) = 0$, it should be an analytical solution to $(BVP)$.
 \hfill $\Box$\\

Generally speaking, the canonical dual algebraic equation  (\ref{eq-cde} ) is nonlinear  which allows multiple solutions
for nonconvex problems. In  the following sections, we shall present some applications.

\section{Application  to Convex Variational Problem}
\begin{figure}[h]
\setlength{\unitlength}{.4cm}
\begin{picture}(-1,5)
\put(3,-6){{\large \epsfig{file=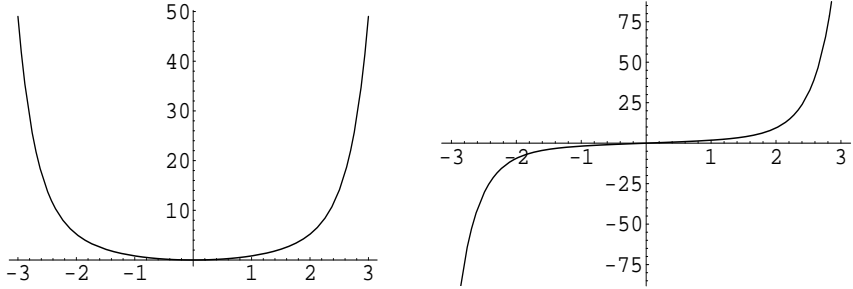,height=4cm,width = 12cm}}}
 \end{picture}\vspace{2.5cm}\\
\caption{Graphs of $\hat{W}(\gamma)$ (left) and its derivative (right) ($\mu=1.0, \; \nu=0.5$)
}\label{admif1}
\end{figure}

First, we assume the stored energy $\hat{W}(\bgamma)$ is a convex function of the type (see \cite{gao-ogden-zamp}):
\begin{equation}
\hat{W}(\bgamma)=\frac{_1}{^2}\mu |\bgamma|^2 + \nu \left(\exp (\half
|\bgamma|^2 ) - 1 \right),
\end{equation}
where $\mu > 0$ and $\nu > 0$ are material constants.
In this case, the constitutive equation (\ref{shear-stress}) can be written as the
following form
\begin{equation}
\btau(\bgamma)= \nabla  \hat{W} (\bgamma)=\mu\bgamma+ \nu  \bgamma \exp
(\half  |\bgamma|^2 ),\label{eq-sg}
\end{equation}
which can  be used to model a large class of materials, especially
 bio-materials (cf. \cite{hol-ogden}).  The  associated potential variational problem is
\eb
(\calP_1): \;\; \min_{\uu \in \calX_a}  \left\{ \Pi(\uu) = \int_\Oo \left[
\frac{_1}{^2}\mu |\nabla \uu |^2 + \nu \left(\exp (\half
|\nabla \uu |^2 ) - 1 \right) \right] \dO    - \int_\Gt \bt \uu \dG  \right\} .
\ee
 This problem also appears in  the construction of optimal Lipschitz extensions of given
boundary data,  the Monge - Kantorovich optimal mass transfer problem, and  a
form of weak KAM theory for Hamiltonian dynamics (see \cite{barron}), etc.
  Although the
energy function is convex   and the
constitutive relation is monotone (see Fig. \ref{admif1}),
the complementary energy $\hat{W}^*(\bsig )$ can not be
obtained by the Legendre transformation since
the inverse relation of $\bsig(\bgamma)$ is analytically
impossible.

  By using the geometrical measure
 $\xi = \Lambda(\uu) =  \half |\nabla \uu |^2  $, the   canonical energy function
 can be defined by
 \eb
\VV(\xi) = \mu \xi + \nu \left(\exp (\xi) - 1 \right),
\ee
and   $\hat{W}(\bgamma) = \VV(\xi(\bgamma))$.
Clearly, the  canonical strain energy $\VV(\xi)$
is well-defined on the domain
$\calE_a$.
Thus, the constitutive law (\ref{eq-sg}) can be written in the following simple form:
\eb
\zeta =   \VV'(\xi) = \mu + \nu \exp(\xi) ,
\ee
 which is uniquely defined on
the domain
\eb
\calE^*_a= \{ \zeta \in {\cal C}[\Oo; \real]  |\; \zeta (\bx) \ge \mu + \nu \;\; \forall \bx \in \Oo \} .
\ee
 Therefore, the complementary energy
$\VV^*(\zeta):\calE^*_a \rightarrow \real$ can be obtained easily as
\eb
\VV^*(\zeta) = \sta \{ \xi \zeta - \VV(\xi) | \;\; \xi \in \calE_a\} =
(\zeta - \mu) \left( \log\left( \frac{\zeta - \mu}{\nu} \right) - 1 \right) + \nu.
\ee
 Clearly, the canonical duality relations
 (\ref{eq-cdr}) hold on $\calE_a \times \calE^*_a$.

By the fact that  on $\calE^*_a$,  we have
$\calS_a = \calE^*_a = \calS_a^+$ and
the canonical stress $\zeta(\bx)  \ge \mu + \nu >  0 \;\; \forall \bx \in \Oo$, the total complementary energy
$\Xi(\uu, \zeta)$ (or $\Xi_\btau(\uu, \zeta)$) is convex in $\uu \in \calX_a$ and
 the pure complementary variational problem $(\calP^d)$ for this convex problem can be written in the following
\eb
(\calP^d_1): \;\;\; \max_{\zeta \in \calS_a}  \left\{ \Pi^d_1(\zeta)
=    - \int_\Oo \left( \;\; \frac{|\bsig|^2   }{2   \zeta}   +
 (\zeta - \mu) \left( \log\left( \frac{\zeta - \mu}{\nu} \right) - 1 \right) + \nu\right)   \dO
  \right\}.
 \ee
 \begin{thm}
  For a given statically admissible  stress  $\btau \in \calT_a$,
   the canonical dual problem $(\calP^d_1)$
  has a unique solution $\bar{\zeta}(\bx)\ge \mu + \nu $.
  If $\btau$ is a potential field and $\btau \times (\nabla \bar{\zeta} ) = 0\;$ on $ \Oo$, then
        the function
  \eb
  \bar{\uu}(\bx) = \int_{\bx_0}^\bx  ( \bar{\zeta} (\bx))^{-1}  \btau \cdot  {\rm d} \bx \;\; \forall \bx_0 \in \Gu \label{eq-solu2}
  \ee
 is a unique solution of  $(\calP_1)$
  and
  \eb
  \Pi_1(\bar{\uu}) = \min_{ \uu \in \calX_a}  \Pi_1(\uu) = \max_{\zeta \in \calS_a } \Pi_1^d(\zeta)
  = \Pi_1^d(\bar{\zeta}).
  \ee
  \end{thm}

\begin{figure}[h]
\setlength{\unitlength}{.4cm}
\begin{picture}(-1,5)
\put(6,-6){{\large \epsfig{file=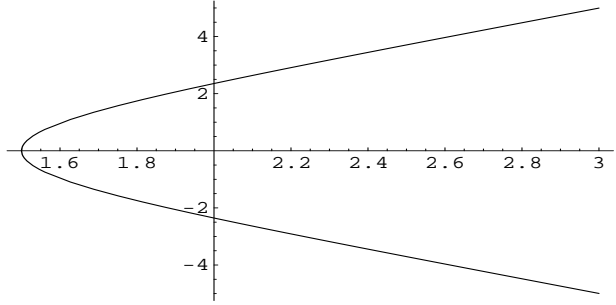,height=4cm,width = 8cm}}}
 \end{picture}\vspace{2.5cm}\\
\caption{Dual algebraic curve $h(\zeta)$ ($\mu = 1, \;\; \nu = 0.5$)
}\label{dac1}
\end{figure}

  {\bf Proof}.
The criticality condition $\delta \Pi_1^d(\zeta) = 0 $
leads to the dual algebraic equation
\eb
 2  \zeta^2 \log \left( \frac{\zeta - \mu}{\nu} \right) =
     |\btau|^2  .\label{eq-dab}
 \ee
Let $h^2(\zeta )  =2  \zeta^2 \log \left( \frac{\zeta - \mu}{\nu} \right)$.
Then the graph of $h(\zeta) = \pm \zeta \sqrt{2 \log((\zeta - \mu)/\nu)}$ is the so-called dual algebraic curve
(see Fig. \ref{dac1}).
Clearly, for any given $|\btau|$, the canonical dual algebraic equation (\ref{eq-dab}) has a unique solution
$\barzeta \in \calS_a$.
  Since   the canonical complementary energy density
  \[
  \psi(\zeta) = \frac{ |\btau|^2   }{2   \zeta} +
 (\zeta - \mu) \left( \log\left( \frac{\zeta - \mu}{\nu} \right) - 1 \right) + \nu
 \]
 is a strictly convex function of  $\zeta$ on the dual feasible space
 $\calS_a$, for a given shear stress  $ \btau $,
 the   pure complementary energy $\Pi_1^d(\zeta)$ is strictly concave
 on $\calS_a = \calS_a^+$. Therefore, the solution $\barzeta$ of the canonical dual
  algebraic equation
  (\ref{eq-dab}) should  a unique global maximizer $\bar{\zeta }  \ge \mu+ \nu$.
 By Theorem \ref{thm-anas}  we know that if the condition $\btau \times (\nabla \barzeta) = 0$, the
 analytic solution $\baruu$ can be determined by (\ref{eq-solu1}).
  \hfill $\Box$ \\

 \section{Application  to Nonconvex Power-Law Material Model}
 In this section the stored strain energy is assume to be a polynomial function  of the
 shear strain $|\bgamma|$:
 \eb
 \hat{W}(\bgamma) = \frac{\mu}{2 b}
 \left[ \left( 1 + \frac{b}{p} ( |\bgamma|^2  -  \epsilon)  \right)^p  - 1 \right], \label{eq-wnconv}
 \ee
 where $\mu   > 0 $ is the infinitesimal shear modulus,
  $p , b> 0$ are  material parameters, and $\epsilon\in \real$  is a given (internal) parameter, which can be viewed as,
  for examples,
  residue strain \cite{gao-mecc99}, dislocation \cite{gao-mms04},  random defects \cite{gao-li-v},
  or input control in functioning materials (see \cite{gao-russell}).
    If $\epsilon = 0$, this is the power-law material model introduced by Knowles in 1977
  \cite{knowles77}, and   in this case,  the energy function
 possesses  the following properties:
 \[
 \hat{W}(0) =   0, \;\; \nabla \hat{W} (0) = 0, \;\; \nabla^2 \hat{W} (0) =   \mu  {\bf I}  \succ 0,
 \]
which are necessary for $ \hat{W} (\bgamma)$ to be a stored energy.
The associated stress in simple shear is
\eb
\btau = \mu \left( 1 + \frac{b}{p} (|\bgamma|^2 - \epsilon ) \right)^{p-1} \bgamma.
\ee

\begin{figure}[h]
\setlength{\unitlength}{.4cm}
\begin{picture}(-1,5)
\put(3,-3){{\large \epsfig{file=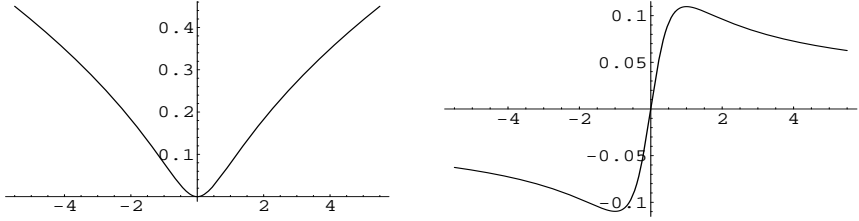,height=3cm,width = 12cm}}}
\put(13,-4.5){(a) $p=0.25, \;\; \epsilon = 0$. }
\put(3,-13.5){{\large \epsfig{file=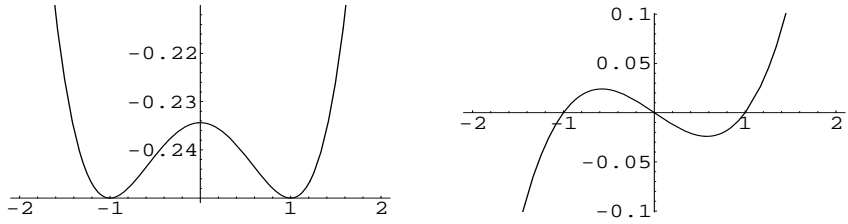,height=3cm,width = 12cm}}}
\put(13,-15){(b) $p=2, \;\; \epsilon =  2.5$. }
\put(3,-24){{\large \epsfig{file=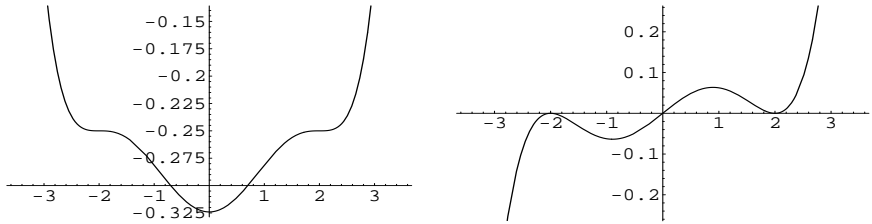,height=3cm,width = 12cm}}}
\put(13,-25.5){(c) $p=3 , \;\; \epsilon = 5$. }
 \end{picture}\vspace{9.5cm}\\
\caption{Nonconvex energy $\hat{W}(\gamma)$ (left) and non-monotone constitutive law $\hat{W}'(\gamma)$
(right), $  \mu = 0.5 ,  \; b=1.$
}\label{fig-admifig2}
\end{figure}

 The power-law material hardens or softens in shear according to whether $p > 1$ or $p < 1$.
 Graphs of this material model are shown in
 Figure  \ref{fig-admifig2}.
 Particularly, when $p = \half, \; \epsilon = 0$, the partial differential equation
 $\nabla \cdot  \nabla \hat{W} (\nabla \uu) = 0$ becomes
 \eb
 (1+ 2 b \uu^2_{,2}) \uu_{,11}  - 4 b u_{,1}\uu_{,2} \uu_{,12} + (1+ 2b \uu^2_{,1})\uu_{,22} = 0,
 \ee
 which, on re-scaling $\uu$ (or by letting $2b = 1$), is the celebrated minimal surface equation
 \eb
 (1+  \uu^2_{,2}) \uu_{,11}  - 2  u_{,1}\uu_{,2} \uu_{,12} + (1+   \uu^2_{,1})\uu_{,22} = 0 .
 \ee
It  also governs the flow of a K\'{a}rm\'{a}n-Tsien gas (see \cite{horgan}).

 It is easy to prove that for $p \ge \half$, the stored energy $\hat{W}(\bgamma)$ is convex (see Section 6.5.3, \cite{gao-dual00}). In this case, the $(BVP)$ is elliptic (see \cite{knowles77}).
 However, if     $p < \half$,    the  constitutive law $\btau  = \nabla \hat{W} (\bgamma)$
 is not monotone even if $\epsilon = 0$    (see Fig.  \ref{fig-admifig2} (a)).
 Although it can  be considered for modeling softening phenomenon,
   this case is not physically  allowed since $p < \half $   violates the constitutive  law.
  Mathematical explanation for this case can be given by the canonical duality theory (see below).

   When $p = 1$  in (\ref{eq-wnconv}), the stored energy function $\hat{W}(\bgamma)$
   is linear, which recovers the neo-Hookean material.

For $p > 1$, the stored energy function $\hat{W}(\bgamma)$   is convex if $\epsilon \le  0$,
nonconvex for $\epsilon > 0$.
Particularly, if $p = 2$ and  $\epsilon > p/b$, this nonconvex function    $\hat{W}(\bgamma)$    is the so-called
 {\em double-well energy} in mathematical physics  (see Fig.  \ref{fig-admifig2} (b)), which appears frequently  in phase transitions of solids, Landau-Ginzburg model in super-conductivity
 \cite{gao-mms04},
 post-buckling of large deformed beam \cite{cai-gao-qin}, as well as in quantum mechanics such as Higgs mechanism
   and Yang-Mills fields etc.
    For $p > 2, $ the stored energy and constitutive law are shown in Fig.  \ref{fig-admifig2} (c).

Let $\beta = \mu b^{p-1}/p^p , \;\; \alpha = \epsilon - p/b$, and $\beta_o = \mu/(2 b)$.
The  minimal potential energy principle leads to the following nonlinear variational extremum problem:
\eb
(\calP_2): \;\; \min \left\{ \Pi_2(\uu) = \int_\Oo
 \left[  \half \beta  \left(  |\nabla \uu |^2  -  \alpha \right)^p  -  \beta_o \right] \dO - \int_\Gt \bt \uu \dG \;\; | \; \uu \in \calX_a \right\}.
 \ee
For certain given parameter $\alpha > 0$, this variational problem
 is nonconvex and the  corresponding boundary-value  problem $(BVP)$ is not equivalent to
 $(\calP_2)$ since the solution of $(BVP)$ may not be the global minimizer of $(\calP_2)$.
Due to the lack of global optimality condition, traditional direct methods for solving the
$(\calP_2)$ are very difficult.

By using the canonical strain measure
  $\Lambda(\uu) =   |\nabla \uu|^2$, the canonical function for this power-law model
 can be defined by
 \eb
 \VV(\xi) =   \half \beta (\xi - \alpha)^p - \beta_o ,
 \ee
which is defined  on the closed convex domain
\[
\calE_a = \{  \xi \in  {\cal L}^p| \; \xi(\bx) \ge 0, \;\;
 \forall \bx \in \Oo \}.
\]
Clearly, this canonical function is convex for $p \ge 1$,   but nonconvex  for $p < 1$.
In any case, the canonical dual stress  is uniquely obtained by
\[
\zeta =    \nabla \VV(\xi)  =  \half p \beta (\xi - \alpha)^{p-1}  ,
\]
which is well defined on the dual space
\[
\calE^*_a = \{   \zeta \in  {\cal L}^{p/(p-1)}   \;  |
\;\;  \zeta (\bx)  \ge  \half \mu (-  \alpha b/p)^{p-1}      \;\; \forall \bx \in \Oo \;\}.
\]
Thus the complementary energy can be simply obtained by the traditional Legendre transformation
\[
\VV^{*}(\zeta ) = \{ \xi \cdot \zeta - \VV(\xi) | \;\; \zeta = \nabla \VV(\xi) \}
=   \frac{p-1}{p} \cc  \zeta^{p/(p-1)} + \alpha \zeta + \beta_o,
\]
where $\cc =  \left( \frac{2 }{p \beta} \right)^{1/(p-1)} = \left( \frac{2}{\mu}\right)^{1/(p-1)} p/b$.
The corresponding total complementary energy for this nonconvex problem  is
\eb
\Xi(\uu, \zeta) = \int_\Oo \left[   \left( |\nabla \uu|^2 - \alpha \right) \zeta
-  \frac{p-1}{p} \cc \zeta^{p/(p-1)} - \beta_o \right]\dO - \int_\Gt \bt \uu \dG.
\ee

Therefore, for a given $\btau \in \calT_a$, let
\eb
\calS_a^+ =  \{   \zeta  \in  \calE^*_a   \;  |
 \;\; \zeta (\bx) \ge 0 \;\; \forall \bx \in \Oo, \;\; |\btau |^2/\zeta \in { \cal L}\;\}.
\ee
The canonical dual problem in this nonconvex case is
 \eb
 (\calP^d_2): \;\;
 \max \left\{ \Pi^d_2(\zeta) = - \int_\Oo \left[ \frac{ \;\; |\btau|^2}{ 4  \zeta } +   \frac{p-1}{p} \cc  \zeta^{p/(p-1)}  + \alpha \zeta  + \beta_o \right]\dO \;
 |\; \zeta \in \calS_a \right\}
 \ee

The criticality condition $\delta \Pi_2^d(\zeta) = 0 $
leads to the dual algebraic equation:
\eb
4    \zeta^2  \left ( \cc   \zeta^{1/(p-1)}     +  \alpha   \right) =   |\btau|^2. \label{eq-dab2}
 \ee
  The solutions of this algebraic equation depends mainly on the material parameter $p > 0$.
  It can be easily checked by using MATHEMATICA that if $p < \half $, this equation has no real root.
  For $p= \half$, the equation (\ref{eq-dab2}) has real roots only under the condition
  $|\btau|^2 \le \frac{\mu^2}{2b}$. Particularly, for minimal surface-type problems
   where $ \mu = 1 $ and $2b = 1$, the condition $|\btau|^2 \le 1$ verifies the result
   presented in \cite{gao-dual00} (Section 6.5.3).  Canonical duality theory for solving
   minimal surface type problems have been studied in
   \cite{gao-yang-95}.
   In this paper, we are interested in $p > 1 $ with positive internal parameter $\alpha > 0$
   such that the stored energy is nonconvex, which can be used to model  more interesting phenomena.

   Particularly, for $p =2$,  the canonical dual algebraic  (\ref{eq-dab2}) is  cubic
\eb
4    \zeta^2  \left ( \cc   \zeta      +  \alpha   \right) =   |\btau|^2 \label{eq-dab3}
 \ee
    which can be solved analytically to have  three solutions:
 \begin{eqnarray}
 \bar{\zeta}_1 &=& - \frac{   \alpha  }{3 c  } +
  \frac{   2^{4/3} \alpha^2 }
  {3 c    \psi(\tau)} +
  \frac{\psi(\tau)}{ 3 (2)^{4/3}  c }  \label{eq-zeta1}\\
\bar{\zeta}_2 &=& -  \frac{   \alpha  }{3 c  }  -
  \frac{ 2^{1/3} \alpha^2  (1 - i \sqrt{3} ) }
  {3 c  \psi(\tau)}
  -  \frac{(1+ i \sqrt{3} )  \psi(\tau)}{12( 2^{1/3} )c },\\
\bar{\zeta}_3 &=&-   \frac{   \alpha  }{3 c  }   -
  \frac{ 2^{1/3} \alpha^2   (1+ i \sqrt{3} ) }
  {3 c \psi(\tau)}
  -  \frac{(1- i \sqrt{3} )  \psi(\tau)}{12 ( 2^{1/3}) c
  },  \label{eq-zeta3}
  \end{eqnarray}
  where $\tau = |\btau|$, and
  \[
  \psi(\tau ) =   \left( -16 \alpha^3 + 27 c^2  \tau^2
   + 3 \sqrt{3 } \;  \tau \; \sqrt{  -32 \alpha^3 c^2 +  27 c^4  \tau^2}  \right)^{1/3}.
  \]

  Similar to the general results proposed in \cite{gao-ima98,gao-na00,gao-ogden-zamp},
 we have the following theorems.

 \begin{thm}[Criteria for Multiple Solutions ]\label{thm-bifur}
 For a given parameter $\alpha \in \real$ and  the  material constants  $p = 2$,  $\mu, \;  b  > 0$ such that $c = 4/(\mu b )> 0$,
 let
 \[
  \eta = \frac{ 16 \alpha^3}{27 c^2}.
  \]
  If $\eta \le  0$,  the $(BVP)$ has a unique solution
 in the whole domain $\Oo$.

 If $ \eta > 0$, then the $(BVP)$ could have multi-solutions
 in $\Oo$. In this case, if $\btau$ is a given shear stress and
  \[
  |\btau(\bx) |^2  >  \eta   \;\; \forall \bx \in \Oo,
 \]
   the dual algebraic equation
 (\ref{eq-dab3}) has a unique real root $\bar{\zeta}(\bx) >  0$.
  If $
 |\btau|^2  <  \eta, $   the dual algebraic equation
 (\ref{eq-dab2}) has three real roots    in the order of
 \eb
\bar{\zeta}_1(\bx) \ge 0 \ge \bar{\zeta}_2(\bx) \ge  \bar{\zeta}_3 (\bx).  \label{rootsod}
 \ee
 \end{thm}

\begin{figure}[h]
\setlength{\unitlength}{.4cm}
\begin{picture}(-1,12)
\put(1,-6){{\large \epsfig{file=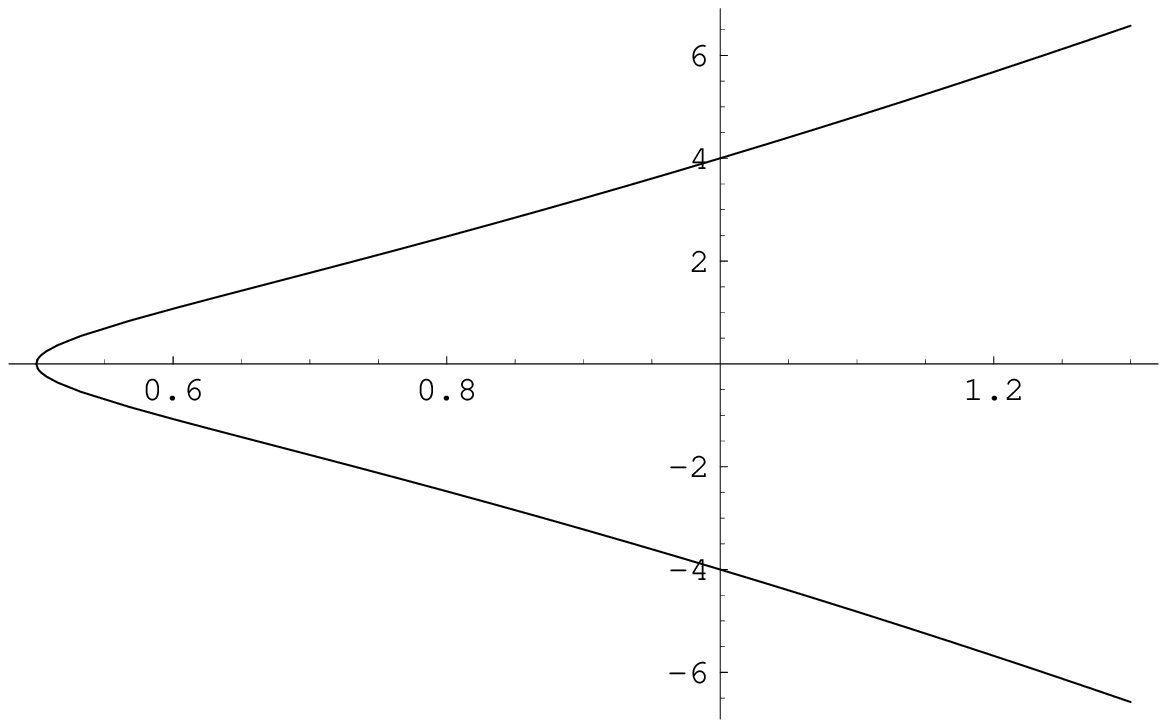,height=7cm,width =5cm}}}
\put(18,-6){{\large \epsfig{file=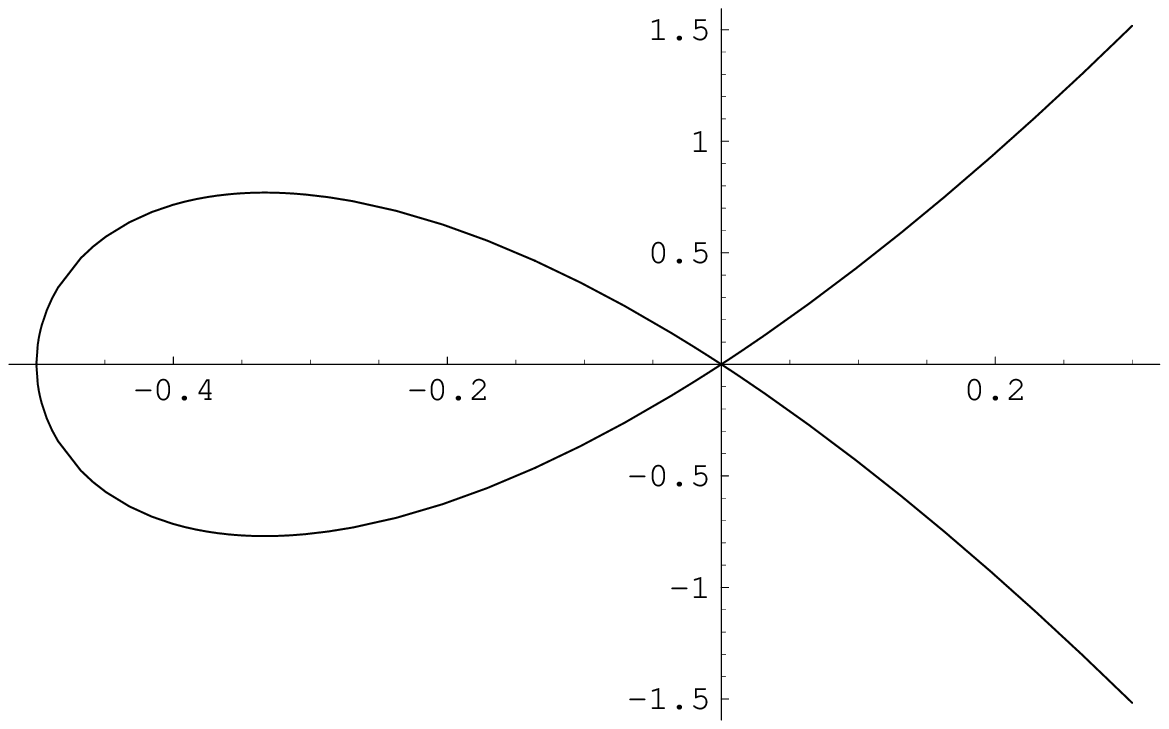,height=7cm,width =6cm}}}
\put(14,2){$ \zeta$}
\put(3.5,11){$ h(\zeta)$}
\put(18,9){\line(1,0){15}}
\put(34,9){$|\btau|^2 >  \eta $}
\put(18,7){\line(1,0){15}}
\put(34,7){$|\btau|^2 =\eta $}
\put(18,5){\line(1,0){15}}
\put(34,5){$|\btau|^2< \eta $}
\put(5,-8){ (a)   $\eta  <  0$ ($\alpha=-4 $). }
\put(22,-8){ (b)  $\eta  > 0$ ($\alpha=4 $ )}
 \end{picture}\vspace{3cm}\\
\caption{Dual algebraic curve   $ h(\zeta) =\pm 2 \zeta  \sqrt{ \cc  \zeta     + \alpha   }  \;
\;( \mu = 1, \;\; b=0.5) $
}\label{fig-admifig3}
\end{figure}

 {\bf Proof}. Similar to the proof of the Corollary 1 in \cite{gao-na00}, we let
  $h^2(\zeta)= 4  \zeta^2 ( \cc \zeta     + \alpha   ) $
  be the left hand side function in
 the dual algebraic equation (\ref{eq-dab3}). By solving $h'(
\zeta_c) = 0$ we known that $h(\zeta)$ has a local maximum
 $h^2_{\max}(\zeta_c) = \eta$ at  $\zeta_c= -\frac{2 \alpha }{3 c}$.
  From the graphs of the dual algebraic curve
  $ h(\zeta) = \pm  2 \zeta   \sqrt{  c  \zeta     + \alpha   }$
   given in Fig. \ref{fig-admifig3}a  we can see that if $\eta < 0$,
   the dual algebraic equation (\ref{eq-dab3}) has a unique solution
   for any given $\btau$.
 However, if    $\eta > 0$,
   the dual algebraic equation (\ref{eq-dab3}) may have at most
   three real  solutions in the order of (\ref{rootsod})
  depending on  $\btau(\bx)  , \;\; \bx \in \Oo$
  (see Fig. \ref{fig-admifig3}b). \hfill $\Box$

 \begin{thm} [Global \& Local Extrema, Uniqueness, and Smoothness]\label{thm-extreme}
Suppose  for a given external force  $\bt (\bx) $ on $\Gt$   that
  $\btau\in \calT_a$ is a statically admissible shear force field.
If $\btau(\bx)  $ is not identical zero over the domain $\Oo$,
 the canonical dual problem $(\calP_2^d)$   has at most three
 solutions  $\bar{\zeta}_i(\bx) \; (i=1,2,3)  $  defined analytically by (\ref{eq-zeta1}-\ref{eq-zeta3}),
 and
 \eb
 \bar{\uu}_i = \int_{\bx_0}^{\bx} (2 \bar{\zeta}_i(\bx))^{-1} \btau \cdot \mbox{d} \bx , \;\; i = 1,2,3 \label{analy3}
 \ee
 are the critical solutions to $(\calP_2)$.
Particularly,  $\barzeta_1(\bx)$
  is a   global maximizer of $\Pi_2^d$ over $  \calS^+_a$, the associated  $ \bar{\uu}_1$ is a global minimizer of $\Pi_2$  on $\calX_a$, and
  \eb
  \Pi_2(\bar{\uu}_1) = \min_{\uu \in \calX_a} \Pi_2(\uu) = \max_{\zeta \in \calS^+_a} \Pi_2^d(\zeta)
  = \Pi_2^d(\bar{\zeta}_1). \label{eq-tri21}
  \ee

 If
$ |\btau (\bx)|^2 < \eta   \;\; \forall \bx \in\Oo,$ then  $\barzeta_3(\bx)$ and the associated  $ \bar{\uu}_3$  are
    local  maximizers of $\Pi_2^d$  and  $\Pi_2$, respectively,   and
  \eb
  \Pi_2(\bar{\uu}_3) = \max_{\uu \in \calX_3} \Pi_2(\uu) = \max_{ -\alp \beta <\zeta < -2\alp \beta/3 } \Pi_2^d(\zeta)
  = \Pi_2^d(\bar{\zeta}_3) ;\label{eq-tri23}
  \ee
   If
$ |\btau (\bx)|^2 < \eta   \;\; \forall \bx \in\Oo \subset \real,$ then  $\barzeta_2(\bx)$ and the associated  $ \bar{\uu}_2$  are
    local  minimizers of $\Pi_2^d$  and  $\Pi_2$, respectively,   and
  \eb
  \Pi_2(\bar{\uu}_2) = \min_{\uu \in \calX_2} \Pi_2(\uu) = \min_{ -2\alp \beta/3<\zeta < 0 \;\;  } \Pi_2^d(\zeta)
  = \Pi_2^d(\bar{\zeta}_2) ,\label{eq-tri22}
  \ee
  where $ \calX_2$ and $\calX_3$ are neighborhoods of  $\bar{\uu}_2$ and $\bar{\uu}_3$, respectively.

   If $ |\btau (\bx)|^2 >  \eta > 0  \;\; \forall \bx \in\Oo,$
 the canonical dual problem $(\calP_2^d)$  has a unique solution $\barzeta_1(\bx)$ over $\Oo$ and
  the primal solution $\bar{\uu}_1$ is a unique smooth global minimizer of $\Pi_2(\uu)$ over $\calX_a$.
  \end{thm}

  {\bf Proof}. If  $\zeta \in \calS^+_a$, the total complementary energy $\Xi(\xi, \zeta)$ is a saddle functional.
  The proof of the canonical min-max  duality (\ref{eq-tri21})  follows directly from Gao and Strang's work \cite{gao-strang89a} and the canonical duality theory
  \cite{gao-dual00}. If  $\zeta < 0 $, the total complementary energy $\Xi(\xi, \zeta)$ is concave in both $\xi$ and $\zeta$.
  The double-max duality (\ref{eq-tri23}) is simply due to the fact
  \[
\max_{\uu \in \calX_a} \Pi_2(\uu) = \max_{\uu } \max_{\zeta  } \Xi_\btau(\uu,\zeta) =
 \max_{\zeta  } \max_{\uu }\Xi_\btau(\uu,\zeta) = \max_{\zeta } \Pi^d_2(\zeta).
 \]

Note that  the double-min duality (\ref{eq-tri22}) holds only for $\Oo\subset \real$. Therefore, by considering   $ \nabla \uu = \uu' = \gamma$
and the canonical transformation  $ \hat{W}(\gamma) = \VV(\xi(\gamma))$, we have
\eb
\nabla^2  \hat{W}(\gamma)  = 2 \nabla \VV(\xi)  + 4 \gamma^2 \nabla^2 \VV (\xi) = 2 \zeta  + \beta (\tau /\zeta)^2  \label{eq-dw2}
\ee
which is positive for any   $\zeta > 0$. Therefore, $\uu_1$ is a global minimizer of $\Pi_2$.
By Theorem \ref{thm-bifur}, it is easy to verify that
\eb
\nabla^2  \hat{W}(\gamma)  \left\{ \begin{array}{ll}
> 0 & \mbox{ if  }  \zeta  \in (- 2 \alp\beta /3 ,  0 ),  \\
= 0 & \mbox{ if  } \tau^2  = \eta, \; \zeta = \zeta_c = - 2 \alp\beta /3, \\
< 0 & \mbox{ if }    \zeta \in ( - \alp\beta , - 2 \alp\beta /3 ) .
\end{array} \right. \label{eq-3w2}
\ee
This shows that $\hat{W}(\gamma) $ is  locally convex at  $\gamma_2 = \tau/(2 \zeta_2)$ and concave
at  $\gamma_3 = \tau/(2 \zeta_3)$.
Therefore $\uu_2$ is a local minimizer, while $\uu_3$ is a local maximizer of $\Pi(\uu)$.

  By the fact that $\tau (\bx) = |\btau(\bx)| \ge 0 \;\; \forall \bx \in \Oo$,
 the force field $\tau(\bx)$ does not cross the Maxwell line, i.e. the $\zeta$-axis in
 Fig. \ref{fig-admifig3} (see Theorem 1 in \cite{gao-ogden-qjmam}).
  Therefore, by Theorem 6 in \cite{gao-ogden-qjmam},
   the global optimal solution should be smooth over the whole domain $\Oo$.
   \hfill $\Box$\\

\begin{remark}[Triality Theory \& Ellipticity Condition]
{\em By Theorems \ref{thm-bifur} and \ref{thm-extreme} we know that if $ |\btau (\bx)|^2 < \eta   \;\; \forall \bx \in\Oo$,
nonconvex problem $(\calP_2)$ has three sets of solutions $\{\bar{\uu}_i(\bx)\} \;(i=1,2,3)$ at each $ \bx \in\Oo$
defined by (\ref{analy3}).
The global minimizer    $\bar{\uu}_1$ is identified by
the canonical min-max duality (\ref{eq-tri21}), which was first proposed  by Gao and Strang in 1989 \cite{gao-strang89a}.
 The (biggest) local maximizer  $\bar{\uu}_3$ is identified by the double-max duality (\ref{eq-tri23}).
Although the canonical dual solution  $ \bar{\zeta}_2$ is a local minimizer of $\Pi^d_2$, the associated $\bar{\uu}_2$ is a local minimizer of
$\Pi_2(\uu)$  governed by the double-min duality (\ref{eq-tri22}) only if the  domain $\Oo$ is a subset of $\real$.
The  {  triality theory} was originally proposed  and proved
for one-dimension problems $\Oo \subset \real$ by Gao in 1996-2000 \cite{gao-amr97,gao-dual00,gao-na00}.
However, in 2003 some counterexamples  were discovered which show that the  double-min duality holds under certain additional constraints
(see Remark 1 in \cite{gao-opt03} and Remark  in \cite{gao-amma03}, page 288).
This open problem has been solved in 2010 first in global optimization, i.e. the ``certain additional constraints" are simply
that the primal and dual problems should have the same dimension in order to have strong triality theory.
Otherwise, the double-min duality holds weakly in a subspace \cite{gao-wu-jimo,gao-wu-jogo}.

Ellipticity condition has been emphasized  to play  a fundamental role for  the existence of solutions in nonlinear elasticity
\cite{fos-ser}. However, this is only for convex problems.
From Theorem \ref{thm-extreme} we know that a nonconvex finite deformation problem could have multiple critical solutions
at each material point and  the associated Euler equation may not be elliptic at all.
Therefore, the  triality theory reveals an important fact in nonconvex analysis and nonlinear elasticity, i.e.
the Legendre-Hadamard  condition does not guarantee uniqueness of solutions and the equilibrium equation may not be elliptic even if
the L.H. holds at certain local solutions.
}
 \end{remark}

\section{Objectivity, Canonical Duality, and Gap Function}
The main goal of this section is to discuss some concepts in   canonical duality theory and  their important roles in
finite elasticity and nonconvex analysis.
Standard notations in 3-D nonlinear elasticity are  adopted.

For a general finite deformation problem  $\bchi : \Oo \subset \real^3 \rightarrow \omega \subset \real^3$,
 the stored energy function $W(\bF)$
is usually a nonconvex function of the deformation gradient $\bF = \nabla \bchi   \in
\mathbb{M}_+^3 = \{ \bF= \{ F^i_\alp\} \in \real^{3\times 3} | \;\; \det F > 0  \}$.
Thus, the boundary-value problem
\eb
(BVP): \;\; \left\{
\begin{array}{l}
A(\bchi) = \nabla \cdot \partial_{\bF} W(\bF(\bchi)) = 0 \;\; \mbox{ in } \Oo,\\
\bn \cdot \partial_{\bF} W(\bF( \bchi) ) = {\bf t} \;\; \mbox{ on } \Gamma_t, \;\; \bchi = \bchi_0 \;\; \mbox{ on } \Gamma_\chi
\end{array}\right.
 \label{eq-ebp}
 \ee
 may have multiple solutions at each material point $\bx \in \Oo$.
According to Dacorogna \cite{daco},  the following statements,  essentially due to   Morrey \cite{morrey},  are well-known:
\eb
\mbox{(I)}  \; \; W(\bF) \mbox{ is  convex } \Rightarrow   \mbox{ poly-convex } \; \Rightarrow
 \; \mbox{ quasi-convex } \; \Rightarrow   \mbox{ rank-one convex}.
\ee

(II) If  $\Oo \subset \real$ or $\omega \subset \real$,
all these notions are equivalent.

(III)  If $W \in {\cal C}^2(\mathbb{M}_+^3)$, then the
rank-one convexity is equivalent to the Legendre-Hadamard (L.H.) condition:
\eb
\sum_{i,j = 1}^3 \sum_{\alp,\beta =1}^3 \frac{\partial^2 W(\bF) }{\partial F^i_\alp \partial F^j_\beta} a_i a_j b^\alp b^\beta \ge 0 \;\; \forall {\bf a} = \{ a_i \}  \in \real^3 , \; \forall {\bf b} = \{ b^\alp \} \in \real^3.
\ee
The   Legendre-Hadamard condition in finite elasticity
is also referred to as the {\em  ellipticity condition}, i.e.,
if the L.H. condition holds, the partial differential operator $A(\bchi)$ in (\ref{eq-ebp})
 is considered to be elliptic.%

However, all these generalized convexities   provide  only necessary conditions
 for local minimal solutions.
 From the triality theory we know that
 the nonconvex variational problem may have multiple solutions at each material point $\bx \in \Oo$.
 The conditions in (\ref{eq-3w2}) show that even if the  Legendre-Hadamard condition  holds  at the solutions $\uu_1(\bx)$ and $\uu_2(\bx)$,
 the   stored energy $\hat{W}(\nabla \uu)$ is not convex at $\bx$ and the
 differential operator $A(\uu)$ is  not monotone, i.e. the statement (II) is not true!
Also, the definition of elliptic operators was originally introduced for
linear partial differential equations that generalize the Laplace equation,
where, the stored energy
 is a convex quadratic function and its  level set   is an ellipse. This definition was generalized for
nonlinear  operators \cite{evans,shubin}.
From the following discussion, we can see that the  stored energy  $W(\bF)$  is not convex  even if the L.H. condition holds.

By the fact that the deformation gradient  $\bF$ is a two-point tensor field, which is not considered as a strain measure.
According to the  {\em axiom  of  objectivity or frame-invariance},
 the following theorem  lays a foundation for the canonical duality theory.
\begin{thm}[Theorem 4.2-1 in \cite{ciarlet-1}]
The stored energy function of a hyper-elastic material  is objective if and only if
\eb
W(\bF) = W(\bQ \bF) \;\; \forall \bF \in \mathbb{M}_+^3, \;\; \forall \bQ \in \mathbb{O}_+^3,
\ee
or equivalently, if and only there exists a function $\VV:\mathbb{S}_>^3 \rightarrow \real$
such that
\eb
W(\bF)  = \VV(\bF^T \bF) \;\; \forall \bF \in \mathbb{M}_+^3,
\ee
where $\mathbb{O}_+^3 = \{ \bQ \in \real^{3\times 3} | \; \bQ^{-1} = \bQ^T, \;\; \det \bQ = 1\}$ is an
 orthogonal  group  and
$\mathbb{S}_>^3 = \{ \bC \in \real^{3\times 3} | \; \bC = \bC^T, \;  \bC \succ 0  \} $.
\end{thm}

The objectivity is a fundamental concept in continuum physics (see \cite{gao-dual00,holz,mars-h,ogden})\footnote{The concept of objectivity has been misused in mathematical optimization  (but mainly in English literature).
  It turns out that
  Gao-Strang's work and the canonical duality-triality theory have been
    challenged  recently by C.    {\em $Z\check{a}$}linescu
  and his co-workers  R. Strugariu,  M. D. Voisei  (see \cite{svz,vz} and references cited therein).
Unfortunately, they  oppositely chose   linear functions as the stored energy $W(\bF)$ and nonlinear functions   as the external energy $F(\bchi)$, they produced many interesting   ``counterexamples" with opposite conclusions. Interested readers are recommended to read \cite{gao-wu-jogo}
  for further discussion.}.
  It was emphasized by P.G. Ciarlet that the  objectivity is an axiom, not an assumption \cite{ciarlet,ciarlet-1}.
  According to the  traditional philosophical
principle of ying-yang duality \cite{gao-96},
 the constitutive relations  in any physical system should be one-to-one in order to obey the
 fundamental law of   nature, i.e. the  Dao  (I-Ching, 2800-2737 BCE).
 This one-to-one constitutive relation  is called the canonical duality.
Therefore, for a  given material, it is reasonable to assume the existence of
 an objective strain measure
 $\bxi = \Lambda(\bchi) \in \calE \subset  \mathbb{S}_>^3 $ and a convex function $\VV:\calE \rightarrow \real$
 such that $W(\bF) = \VV(\bxi)$ and the following canonical duality relations hold
 \eb
 \bxi^* = \nabla \VV(\bxi) \;\; \Leftrightarrow \;\; \bxi = \nabla \VV^*(\bxi^*) \;\;
 \Leftrightarrow \;\; \VV(\bxi) + \VV^*(\bxi^*) = \bxi : \bxi^* .
 \ee
 By the canonical transformation $W(\nabla \bchi) = \VV(\Lambda(\bchi))$, the general
  minimal potential problem in finite deformation theory can be written in the
 following canonical form
 \eb
 \min \left\{ \Pi(\bchi) = \int_\Oo \VV(\Lambda(\bchi)) \dO + F(\bchi) \; | \; \bchi \in {\cal X}_a \right\} , \label{eq-gao-st}
 \ee
which is the   mathematical model studied by Gao and Strang in \cite{gao-strang89a}
for general geometrically nonlinear systems,
 where $F(\bchi)$ is the so-called external energy, which should be a linear functional of $\bchi$
 such that its G\^{a}teaux  derivative obeys the
 Newton third law of action and reaction;  the feasible space is defined by
 \eb
 {\cal X}_a = \{ \bchi \in {\cal C}[\bar{\Oo}; \real^3] | \;
 \nabla \bchi \in \mathbb{M}_+^3, \;\; \bchi = \bchi_0 \mbox{ on } \Gamma_{\chi} \}.
 \ee

 Canonical duality theory has been extensively studied for different objective  measures $\bxi = \Lambda(\bchi)$
 in   continuum mechanics and general complex systems \cite{gao-dual00,gao-amma03,gao-sherali}.
 In order to understand why the complementary  gap function can be used to identify both global and local extrema,
  let us consider the most simple   canonical strain measure
 $\bxi = \half   \bF^T \bF$ such that $\bE = \bxi - \half \bI = \half (\bF^T \bF -\bI)$ is the    well-known Green-St Venant strain  strain.
 Its canonical dual is the second Piola-Kirchhoff stress, denoted by $\bT = \nabla \VV(\bE)$.
 For a given statically admissible stress $\btau \in \calT_a = \{ \btau \in \real^{3\times 3} | \; \nabla \cdot \btau (\bx) = 0 \;\; \forall \bx \in \Oo, \;\; \bn\cdot \btau (\bx)= {\bf t} \; \forall \bx \in \Gamma_t\}$, the pure complementary energy can be formulated as
 \eb
 \Pi^d(\bT) = \int_{\Gamma_\chi} \bn \cdot \btau \cdot \bchi_0 \dG -  \int_\Oo \frac{1}{2} \tr ( \btau \cdot  \bT^{-1} \cdot \btau + \bT ) \dO - \int_\Oo \VV^*(\bT) \dO.
 \ee
 The criticality condition $\delta \Pi^d(\bT) = 0 $ leads to the canonical dual algebraic equation \cite{gao-mrc99}
 \eb
  \bT \cdot ( 2 \nabla \VV^*(\bT) + \bI ) \cdot \bT = \btau^T \btau. \label{eq-cde3}
 \ee
 Clearly, for a given statically admissible stress field $\btau(\bx) \in \calT_a$, this nonlinear tensor equation may have multiple solutions
 $\{\bT_k \}$, and for each of these critical solutions,  the deformation defined by
 \eb
 \bchi_k (\bx) = \int_{\bx_0}^\bx  \btau \bT_k^{-1} \mbox{d} \bx + \bchi_0(\bx_0)
 \ee
along any path from $\bx_0 \in \Gamma_{\chi}$ to $\bx \in \Oo$ is a critical point of $\Pi(\bchi)$ \cite{gao-dual00}.
The vector-valued function
$ \bchi_k (\bx)$ is a solution to the boundary-value problem $(BVP)$ if the compatibility condition
$ \nabla \times (\btau \cdot \bT_k^{-1} ) = 0  $ holds \cite{gao-mrc99}.
 By Gao-Strang's work \cite{gao-strang89a},  $ \bchi_k (\bx)$ is a global minimizer of $\Pi(\bchi)$  if
  the complementary  gap function
 \eb
 G(\bchi, \bT_k) = \int_\Oo \half \tr [(\nabla \bchi)^T \cdot \bT_k \cdot (\nabla \bchi) ]\dO \ge 0 \;\; \forall \bchi \in {\cal X}_a . \label{eq-gap3}
 \ee
Since this gap function is  quadratic in $\bchi$,
  the sufficient condition (\ref{eq-gap3}) holds  if $\bT_k \in \mathbb{S}^3_> $.
 By the triality theory, $ \bchi_k (\bx)$ could be a local minimizer   if $\bT_k \in \mathbb{S}^3_< $.
 To see this, let us consider the general case of the equation (\ref{eq-dw2}), i.e., by chain role for $W(\bF) = \VV(\bE(\bF))$, we have
 \eb
 \frac{\partial^2 W(\bF)}{\partial F^i_\beta \partial F^j_\beta} = \delta^{ij}  T_{\alp\beta} +
 \sum_{\theta, \nu = 1}^3  F^i_\theta  H_{\theta \alp\beta \nu} F^j_\nu,
 \ee
 where ${\bf H} = \{ H_{\theta \alp\beta \nu}\} = \nabla^2 \VV(\bE)$. By the convexity of  the canonical function $\VV(\bE)$,
the Hessian ${\bf H}$ is positive definite.  Therefore, if
  $\bT = \{ T_{\alp\beta} \}\in \mathbb{S}^3_> $,  the
 L.H. condition holds.  By the fact that $\bF = \btau \bT^{-1}$,
  if $\bT  \in \mathbb{S}^3_<  = \{ \bT\in \real^{3\times 3} | \; \bT= \bT^T, \;  \bT  \prec  0  \}  $, the L.H. condition could  also hold depending on the eigenvalues of  $\bT$.
Therefore, if  the boundary-value problem (\ref{eq-ebp}) has multiple solutions $\{ \bchi_k (\bx) \}$
at $\bx \in \Oo$, the operator
 $A(\bchi) $ may not be elliptic at $\bx \in \Oo$ even if the L.H. condition holds  at certain  $\bchi_k (\bx) $.

For St Venant-Kirchhoff material, the strain energy $\VV(\bE)$ is convex (quadratic)
 \eb
 \VV(\bE) = \mu \tr (\bE^2) + \half \lambda (\tr \bE)^2,
 \ee
 where $\mu, \lambda > 0$ are Lam\'{e} constants.
 In this case,  the complementary energy  is
 \eb
 \VV^*(\bT) = \frac{1}{4 \mu} \tr (\bT^2)  - \frac{\lambda}{4 \mu ( 3 \lambda + 2 \mu)} (\tr \bT )^2
 \ee
  and
 the canonical dual algebraic equation (\ref{eq-cde3})
  is a cubic symmetrical tensor equation
 \eb
 \bT^2 + \frac{1}{\mu} \bT^3 - \frac{\lambda}{\mu ( 3 \lambda + 2 \mu)} (\tr \bT ) \bT^2 = \btau^T \btau.
 \ee
 It was shown in \cite{gao-h-ogden} that for a given $\btau (\bx) \neq 0$, this canonical dual algebraic equation has
 a unique positive-definite  solution $\bT_+ \in \mathbb{S}^3_>$, at most nine negative-definite solutions
 $\bT_- \in \mathbb{S}^3_< $, and at most 17 indefinite solutions at each material point $\bx \in \Oo$.
 By the triality theory, $\bT_+ \in \mathbb{S}^3_>$ gives the global minimizer of the total potential $\Pi(\bchi)$;
the smallest $\bT_- \in \mathbb{S}^3_< $ leads to local maximizer, while the
 biggest $\bT_- \in \mathbb{S}^3_< $  could give a local minimizer of $\Pi(\bchi)$.
 Detailed discussion is given in \cite{gao-h-ogden}.

\section{Concluding  Remarks and Open Problems}
Concrete applications of the canonical duality-triality theory have been presented in this paper for solving general anti-plane shear problems in finite elasticity.
Results  show  that the nonconvex variational problem could have multiple
  solutions at each material point $\bx \in \Oo$, the Euler equation is not  elliptic  and
   the Legendre-Hadamard condition is only a local criterion which can't guarantee uniqueness of solutions.
By using the  pure complementary energy principle proposed in \cite{gao-mrc99,gao-mecc99},
the nonlinear partial differential equation in finite elasticity is equivalent to an algebraic (tensor)  equation,
which can be solved, under certain conditions, to obtain a complete set of solutions  in stress
space.  Therefore, an unified analytical solution form is obtained for the nonconvex variational problem.
  The Gao-Strang complementary gap function and the triality theory can be used to identify both global and local extrema.
  By the fact that   the statically admissible
 stress  field $\btau \in \calT_a$ may not be uniquely determined for a given external force ${\bf t}  (\bx) $ on $\Gamma_t$,
  the compatibility condition (\ref{eq-compatible}) should be satisfied in order that this analytical solution
 solves also the mixed boundary-value  problem.
 How to  satisfy this compatibility condition and to identify local minimizer for  3-D problems  are
  still  open questions and  deserve  future study.

\subsection*{Acknowledgements}
The topic of this paper was suggested by Professor David Steigmann at UC-Berkeley.
Important comments and constructive suggestions from three anonymous referees are sincerely acknowledged.
The research   was supported by US Air Force Office of Scientific Research (AFOSR FA9550-10-1-0487).

\end{document}